\documentclass[11pt]{article}

\usepackage{amsmath,amsthm,amsfonts}

\usepackage{amssymb}

\usepackage{epsfig}

\usepackage{rotating}

\usepackage{subfigure}

\newcommand{\p}[2]{\frac{\partial#1}{\partial#2}}

\begin{document}

\title{Numerical computation of solutions of the critical nonlinear Schr\"odinger equation after the singularity}
\author{Panos Stinis \\ 
Department of Mathematics \\
University of Minnesota \\
    Minneapolis, MN 55455} 
\date {}

\maketitle

\begin{abstract}
We present numerical results for the solution of the 1D critical nonlinear Schr\"odinger with periodic boundary conditions and initial data that give rise to a finite time singularity. We construct, through the Mori-Zwanzig formalism, a reduced model which allows us to follow the solution after the formation of the singularity. The computed post-singularity solution exhibits the same characteristics as the post-singularity solutions constructed recently by Terence Tao. 
\end{abstract}


\section*{Introduction}
Nonlinear Schr\"odinger equations with power nonlinearities have been an intense subject of research both analytically and numerically (see e.g. \cite{sulem,tao_book}  and references therein). Depending on the sign and the order of the nonlinear term, the solutions of such equations can exhibit a varied range of behaviors from scattering to solitons to finite time singularities. In the current paper we are interested in the critical Schr\"odinger equation in one spatial dimension. It has been shown that given large enough initial data, the solution of this equation can exhibit finite time blow-up \cite{merle,perelman}. We focus here on the 1D case  because it facilitates the numerical analysis and because the numerical results we present are novel even in this case. However, our techniques can be applied to cases with more spatial dimensions and different degrees of the power nonlinearity. 

We are particularly interested in investigating numerically the behavior of solutions to the equation {\it after} the formation of the singularity. There exists a large body of numerical research on the behavior of the solutions as they approach the blow-up instant \cite{sulem}. Yet, to the best of our knowledge there is no prior numerical work done on what happens to the solution after the singularity and whether it is possible to follow the solution after the singularity has formed. 

The reason that makes the computation of the solution at the singularity (and possibly after) difficult, is that the solution, loses its smoothness at the singularity. It may lose its smoothness for later too (as it happens with shocks) but even the loss of smoothness at one instant is enough to cause severe numerical problems. Loss of smoothness means that there is propagation of activity down to the zero length scale. Thus, the question is how can one represent accurately such a solution since any numerical calculation can only afford to resolve a finite number of length scales (called resolution hereafter). 

We have chosen to address this issue through dimension reduction (see e.g. \cite{givon} for a review) and in particular the Mori-Zwanzig formalism \cite{CHK00}. The main idea behind dimension reduction is to divide the available resolution into resolved and unresolved variables (length scales in our case) and construct a model for the resolved variables. The remaining computational capability is used to transfer activity from the resolved to the unresolved variables. The effect of the unresolved variables on the model equations for the resolved variables is to enhance these equations by terms which account (inevitably in an approximate manner) for the interaction between resolved and unresolved variables. In essence, what one looks for is a way to simplify the dynamics, since a finite resolution does not allow keeping all the dynamics, while at the same time retaining the most important features of the dynamics. As is expected \cite{CS05}, the main difficulty in constructing a reduced model is to estimate the correct rate at which activity is transferred between the resolved and the unresolved variables.  For the nonlinear Schr\"odinger equation we have chosen to use a reduced model stemming from the Mori-Zwanzig formalism called the $t$-model (see \cite{CHK3,bernstein,HS06} for thorough discussions and other applications of this model). 

Recently, Tao \cite{tao} (see also \cite{bourgain,fibich,merle2}) has constructed solutions for the critical nonlinear Schr\"odinger equation in at least four spatial dimensions which hold even after the formation of the singularity. The main feature of these solutions is that they eject a finite amount of mass instantaneously at the singularity instant. These solutions have at most a finite number of mass ejection events (depending on the magnitude of the initial data). Between these events the solution continues to conserve mass. For the special case of spherical symmetry, the amount of mass ejected is at least that of the ground state mass. The solutions of the $t$-model in 1D do eject mass in a very narrow time interval around the singularity instant. The amount of mass ejected is smaller than the ground state mass. It is not known yet rigorously how much mass should be ejected in the 1D case (see Theorem 3 and related comments in \cite{merle2}). 

In addition to mass ejection at the singularity, the solution can lose its regularity for later times \cite{tao,merle2}. Bourgain and Wang \cite{bourgain} have constructed solutions to the 1D critical Schr\"odinger with algebraic blow-up rate which retain their smoothness after the singularity. However, these solutions are unstable to perturbations and have not been observed numerically. There exist also more stable blow-up solutions (with a log-log correction to the blow-up rate) \cite{perelman,merle2} which do lose smoothness after the singularity. The initial condition we have chosen gives rise to a  blow-up solution according to the log-log scenario. This was confirmed by an independent calculation based on a mesh refinement algorithm developed recently by the author \cite{S09}. Thus, we expect a good reduced model to be able to show the ensuing loss of smoothness after the singularity. Even though any numerical calculation has finite resolution, the reduced model should behave in a way that is consistent with the theoretical predictions. Indeed, the solution of the $t$-model after the singularity shows an increasing roughness with increasing resolution. This trend suggests that the solution does indeed lose its regularity after the singularity.  

The paper is organized as follows. Section \ref{schrodinger} presents some generalities about the 1D critical Schr\"odinger equation. In Section \ref{tmodel} we give a very brief presentation of the $t$-model (see \cite{HS06} for more details). Section \ref{numerical} contains the numerical results. In Section \ref{mass_hamiltonian_evolution} we examine the evolution with time of the mass, the Hamiltonian and the $l_2$ norm of the gradient of the solution. In Section \ref{smoothness} we examine the smoothness of the solution after the singularity. Finally, Section \ref{discussion} concludes with a discussion of the results and some directions for future work.


\section{The critical focusing nonlinear Schr\"odinger equation}\label{schrodinger}

The 1D critical focusing Schr\"odinger equation \cite{sulem} is given by 
\begin{equation}
i \p{u}{t} + \Delta u + |u|^4u=0.
\label{schrod} 
\end{equation}
The equation needs to be supplemented by an initial condition $u(x,0)=u_0(x)$ and boundary conditions. We solve (\ref{schrod}) in the interval $[0,2\pi]$ with periodic boundary conditions. If we assume that the solution remains smooth, it is straightforward to show that the solution of \eqref{schrod} conserves the mass $M(t)$ given by 
$$M(t) =\int_{[0,2\pi]} |u(x,t)|^2 dx$$ 
as well as the Hamiltonian $H(t)$ given by 
$$H(t)=\int_{[0,2\pi]}  \biggl[ \frac{1}{2} |\nabla u(x,t)|^2 - \frac{1}{6} |u(x,t)|^6 \biggr] dx, $$
The use of periodic boundary conditions allows us to expand the solution in Fourier series
$$u^{K}(x,t )=\underset{k \in F \cup G}{\sum} u_k(t) e^{ikx},$$
where $F \cup G=[-\frac{K}{2},\frac{K}{2}-1].$ We have written the set of Fourier modes as the union of two sets 
in anticipation of the construction of the reduced model comprising only of the modes in $F=[-\frac{N}{2},\frac{N}{2}-1],$ where $ N < K.$
The equation of motion for the Fourier mode $u_k$ becomes
\begin{equation}
\label{schrodode}
 \frac{d u_k}{dt}= -i k^2 u_k+ i \underset{k_1 \ldots, k_5 \in F \cup G}{\underset{k_1-k_2+k_3-k_4+k_5=k  }{ \sum}} u_{k_1} u^{*}_{k_2} u_{k_3} u^{*}_{k_4} u_{k_5},
\end{equation}
where $u^{*}_k$ denotes the complex conjugate of the Fourier mode $u_k.$ The ODE system \eqref{schrodode} conserves the discrete versions of the mass and the Hamiltonian 
$$ M(t)= \sum_{k \in F \cup G} |u_k(t)|^2   \; \text{and} \;    H(t)= \sum_{k \in F \cup G}  \biggl[ \frac{1}{2} k^2 |u_k(t)|^2 -  \frac{1}{6} |u_k(t)|^6 \biggr].$$

\section{The $t$-model}\label{tmodel}

The solution of \eqref{schrod} can blow-up in finite time depending on the magnitude of the initial condition $u_0(x).$ 
The representation of the solution of \eqref{schrod} by a finite number of Fourier modes breaks down at the blow-up instant since the solution develops activity down to the zero-scale. This presents a major problem for numerical calculations since we can only afford a finite number of Fourier modes. In other words, no matter how large a calculation we can afford, we are bound to run out of resolution at the blow-up instant. However, in some cases, it is possible to construct a model for a reduced set of Fourier modes (called the resolved modes) which remains well-resolved even after the blow-up instant. Such a model needs to be able to eject mass {\it at the correct rate} from the resolved to the unresolved modes. As expected, the main difficulty in constructing the reduced model lies in estimating the correct rate of mass ejection from the resolved to the unresolved modes. In general, the problem of estimating the rate at which activity propagates from the large scales to the small scales of the solution is hard (see \cite{CS05,givon} for extensive discussions and examples).

The ODE system \eqref{schrodode} for the modes $F \cup G$ can be rewritten as 
$$\frac{du(t)}{dt} = R (u(t)),$$
where $u = ( \{u_k\}), \; k \in F \cup G$ is the vector of Fourier coefficients of $u.$ Also, $R (u(t))$ is the vector of right hand sides (RHS) of \eqref{schrodode}. The vector of Fourier coefficients can be written as $ u = (\hat{u}, \tilde{u}),$ where $ \hat{u}$ are the resolved modes (those in $F$) and $\tilde{u}$ the unresolved ones (those in $G$). Similarly, for the RHS vector we have $R(u) = (\hat{R}(u), \tilde{R}(u)).$ 

We need to choose a reduced model for the modes in $F.$ We use a reduced model, known as the $t$-model, which has been shown to follow correctly the behavior of the solution to the inviscid Burgers equation even after the formation of shocks \cite{bernstein,HS06}. It has also been used to investigate the possible finite time blow-up for the 3D Euler equations of fluid mechanics \cite{HS06}. The $t$-model was first derived in the context of statistical irreversible mechanics using the Mori-Zwanzig formalism \cite{CHK3} and was later analyzed in \cite{bernstein,HS06}. It is based on the assumption of the absence of time scale separation between the resolved and unresolved modes. For a mode $u_k'$ in $F$ the model is given by
\begin{multline}\label{schrodode2}
\frac{d}{dt}u_k'=-i k^2 u'_k+ i \underset{k_1 \ldots, k_5 \in F}{\underset{k_1-k_2+k_3-k_4+k_5=k  }{ \sum}} u'_{k_1} u'^{*}_{k_2} u'_{k_3} u'^{*}_{k_4} u'_{k_5} \\
+t \biggr[ 3i\underset{k_1 \in G  ,\, k_2,\ldots,k_5 \in F}{\underset{k_1-k_2+k_3-k_4+k_5=k  }{ \sum}}    R_{k_1}(\hat{u}') u'^{*}_{k_2} u'_{k_3} u'^{*}_{k_4} u'_{k_5}  \\
+2i \underset{k_2 \in G ,\, k_1,k_3,\ldots,k_5 \in F }{\underset{k_1-k_2+k_3-k_4+k_5=k  }{ \sum}}  u'_{k_1}R^{*}_{k_2}(\hat{u}')  u'_{k_3} u'^{*}_{k_4} u'_{k_5} \biggl]. 
\end{multline}
where we have suppressed the dependence of the solution on time to avoid clutter in the formulas. Also, the prime is used to denote the fact that the solution of the reduced model \eqref{schrodode2} can differ from the solution of the system \eqref{schrodode}. Of course, we hope that the reduced model will be able to reproduce the correct behavior for the resolved modes. 

Note that the $t$-model is closed in the resolved modes. The second term on the RHS of (\ref{schrodode2}) is of the same form as the second term in (\ref{schrodode}), except that the term in (\ref{schrodode2}) is defined only for the modes in $F.$ The third and fourth terms in (\ref{schrodode2}) are not present in (\ref{schrodode}). They are of order nine in the Fourier modes and they are effecting the drain of mass out of the modes in $F.$ The name of the $t$-model comes from the explicit dependence on time of the RHS of the system \eqref{schrodode2}. 

The $t$-model has two features which make it attractive as a reduced model for problems without time scale separation over a large range of modes (scales). The first feature is that it is derived directly from the system of ODEs and does not involve adding regularizing terms by hand. This facilitates the proof of its convergence with increasing resolution to the solution of the PDE as long as the solution is smooth. Also, it does not involve any adjustable parameters in contrast to artificially added regularizing terms. The second welcome feature of the $t$-model is that, for systems of ODEs which conserve the $l_2$ norm of the solution (the mass in the nonlinear Schr\"odinger case), the $l_2$ norm for the resolved modes is {\it non-increasing} in time \cite{HS06}. In fact, for $ M_F(t)= \sum_{k \in F } |u'_k(t)|^2$ we have
\begin{equation}\label{mass_dissipation}
 \frac{d M_F(t)}{dt} = -2t\sum_{k \in G } |R_k(\hat{u}'(t))|^2.
 \end{equation}

We should note that the $t$-model has a similar form for the critical nonlinear Schr\"odinger in more than one spatial dimensions. It can also be constructed for the supercritical Schr\"odinger equation (in one or more spatial dimensions). We have applied the $t$-model in these cases and we will present those results in a future publication. 

\section{Numerical results}\label{numerical}

In the numerical experiments we used the (purely imaginary) initial condition 
$$ u_0(x,0)=i A \exp(-(x-\pi)^2).$$
For this initial condition we have $\max |u_0(x)| = A$ at $x=\pi.$ We present results for the case of $A=1.80$ which was found through a mesh refinement algorithm \cite{S09} to lead to a finite-time blow-up (more details on the behavior at blow-up are given later).

Before we present the numerical results we should comment on the choice of the range of the resolved modes (those in $F$) and the unresolved modes (those in $G$). As can be seen from the RHS of the reduced model \eqref{schrodode2}, for the critical Schr\"odinger in one spatial dimension, the $t$-model allows interactions of a mode with modes which have at most {\it five} times larger wavenumber. This means that the optimal division of an available resolution (available number of Fourier modes for a simulation) is to take the resolved range as one fifth of the available range. For example, if we can afford to calculate with, say $K=50$ Fourier modes, we should construct a reduced model for $N=10$ modes and use the other $40$ modes as unresolved. This would mean that $F=[-5,4]$ and $G=[-25,-6] \cup [5,24].$ In the figures, the number $N$ of modes denotes the number of resolved modes. For example, $N=256$ means that the calculation of the $t$-model involves $256*5=1280$ modes. The form of the $t$-model in Fourier space makes it possible to use FFT to calculate the nonlinear sums. If we opt for the division of modes just discussed above, then all the FFTs are dealiased by construction (see \cite{HS06} for more details).  

The largest resolution we have used is $N=512$ which means that the calculation involves 2560 Fourier modes. At first sight this may seem a small resolution for a 1D calculation. However, the situation is more complicated. We have used the Runge-Kutta Fehlberg method of order 4-5 to integrate the equations of the $t$-model with the error control tolerance set to $10^{-10}$ \cite{hairer}. This leads to a stepsize of about $10^{-5}$ before the singularity. As is known, the error control depends on the magnitude of high order temporal derivatives. So, at the singularity, where the solution changes very rapidly in time, the stepsize will be decreased. Indeed, the stepsize falls to around $10^{-7}$ shortly before and after the singularity and then plateaus to about $10^{-6}$ for the remaining calculation. If one wants to use larger resolution, even in 1D, it is advisable to parallelize the algorithm.   

\begin{figure}
\centering
\includegraphics[width=3.5in]{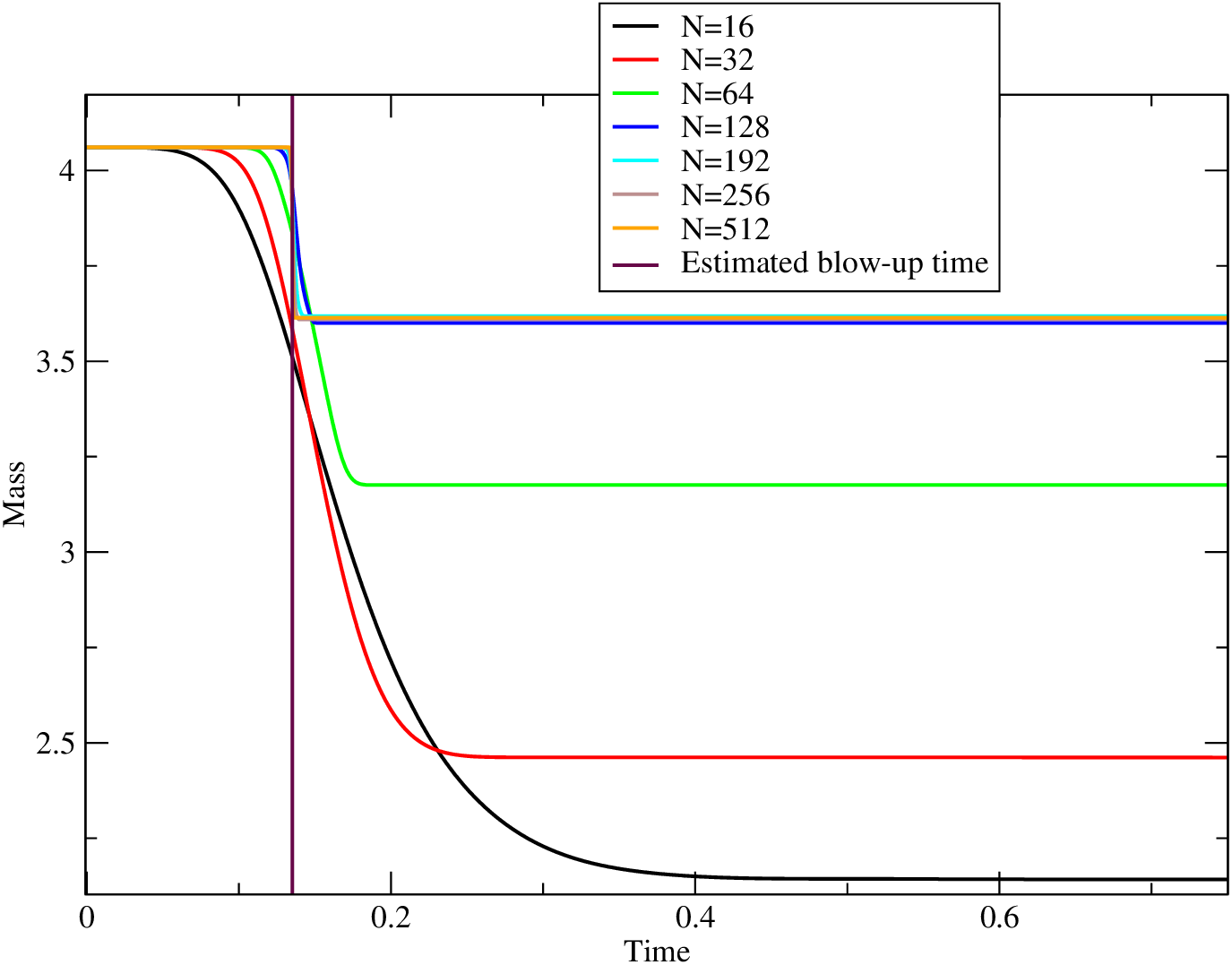}
\caption{Evolution of mass for resolved modes. The vertical line denotes the numerically estimated blow-up instant $T=0.13504$ calculated with a mesh refinement algorithm.}
\label{critical_1}
\end{figure}

\subsection{Mass and Hamiltonian evolution}\label{mass_hamiltonian_evolution}

Figure \ref{critical_1} shows the evolution of the mass for the resolved modes for different resolutions, from $N=16$ to $N=512.$ The behavior of the mass evolution leads us to a few observations. The first observation is that as we increase the resolution the mass of the resolved modes remains constant for a longer time. We know from theory that the mass of the solution of the critical Schr\"odinger equation is conserved for all times that the solution remains smooth. However, the moment the solution blows up, the solution loses its smoothness and there is no reason why the mass should be conserved. In fact, Tao \cite{tao} (see also \cite{bourgain}) has constructed solutions for the critical Schrd\"odinger equation which at the blow-up instant eject mass instantaneously and then continue with a constant (but lower) mass. 

In addition, Tao showed that the number of mass ejection occurrences is finite. While in the current example we have found only one mass ejection occurrence we have conducted numerical experiments with larger initial conditions which give rise to more than one mass ejection occurrences, but always a finite number of them. This is in direct contrast to the behavior of, say, the inviscid Burgers equation where once a shock has been established, the loss of mass (called energy in the relevant literature) is not instantaneous but persists in time and in fact it follows a power law \cite{lax}.

The second observation is that the time of occurrence of the mass ejection as predicted by the $t$-model is in remarkable agreement with the estimated blow-up time based on a mesh refinement algorithm developed earlier by the author \cite{S09}. In other words, the $t$-model kicks in only when needed and does not eject mass unnecessarily. For our initial condition, the algorithm in \cite{S09} shows that the solution blows up following the log-log scenario proven by Perelman \cite{perelman}. In fact,
$$  |\nabla u(x,t)|_{L_2} \sim \biggl( \frac{\ln|\ln(T-t)|}{T-t} \biggr)^\frac{1}{2}    \; \text{for} \; t \rightarrow T^-,$$
where $T$ is the blow-up instant. For $A=1.80,$ the estimated blow-up time $T=0.13504.$ We monitored the mass dissipation rate of the $t$-model given by equation \eqref{mass_dissipation} above and found it to be sharply peaked at $T^*=0.13508,$ which also corroborates the agreement of the $t$-model behavior with the occurrence of a singularity predicted by the mesh refinement algorithm. 

The third observation is that the value of the mass after the singularity, as predicted by the $t$-model,  appears to be converging as we increase the resolution of the reduced model. Of course, one example is not enough to infer the converging properties of the $t$-model especially since we do not know if it converges to the right solution after the singularity.

A related issue is the amount of mass concentrated at the blow-up point and how much mass is ejected at the singularity instant. In \cite{ogawa} it was shown that in 1D and for periodic boundary conditions, the solution concentrates at the blow-up point mass equal to the mass of the ground state on the whole line. The ground state equation is given by \cite{sulem} 
$$ \Delta Q + |Q|^4Q - Q=0.$$
In 1D, its solution on $\mathbb{R}$ is given by 
$$Q(x)=\biggl(  \frac{3}{\cosh^2(2x)}  \biggr)^{\frac{1}{4}}.$$ 
The mass of the ground state  $M_Q =\int_{\mathbb{R}}|Q(x)|^2 dx$ was found to be approximately 2.7412. 

We calculated the mass concentration around the blow-up point which is at $\pi.$ We find that at the instant of the singularity, a mass amount equal to the ground state mass is concentrated in the region $[\pi-0.05,\pi+0.05].$ Of course, since we can only afford a finite resolution, the area in which the ground state mass is concentrated is not a single point but lies in a narrow range around the blow-up point. After the singularity has occurred we find that the mass ejected by the $t$-model (with $N=512$) is about 0.446. The amount of mass ejected is small compared to the ground state mass. Merle and Raphael \cite{merle2} talk about radiative mass ejection for the log-log blow-up scenario, however, at this point we do not know how large this mass ejection should be (see Theorem 3 in \cite{merle2} and the related comments). In \cite{tao}, Tao has shown that under spherical symmetry and in at least four spatial dimensions the amount of mass ejected should be at least the mass of the ground state. However, it is not known yet whether these results apply to lower spatial dimensions too.  

\begin{figure}
\centering
\includegraphics[width=3.5in]{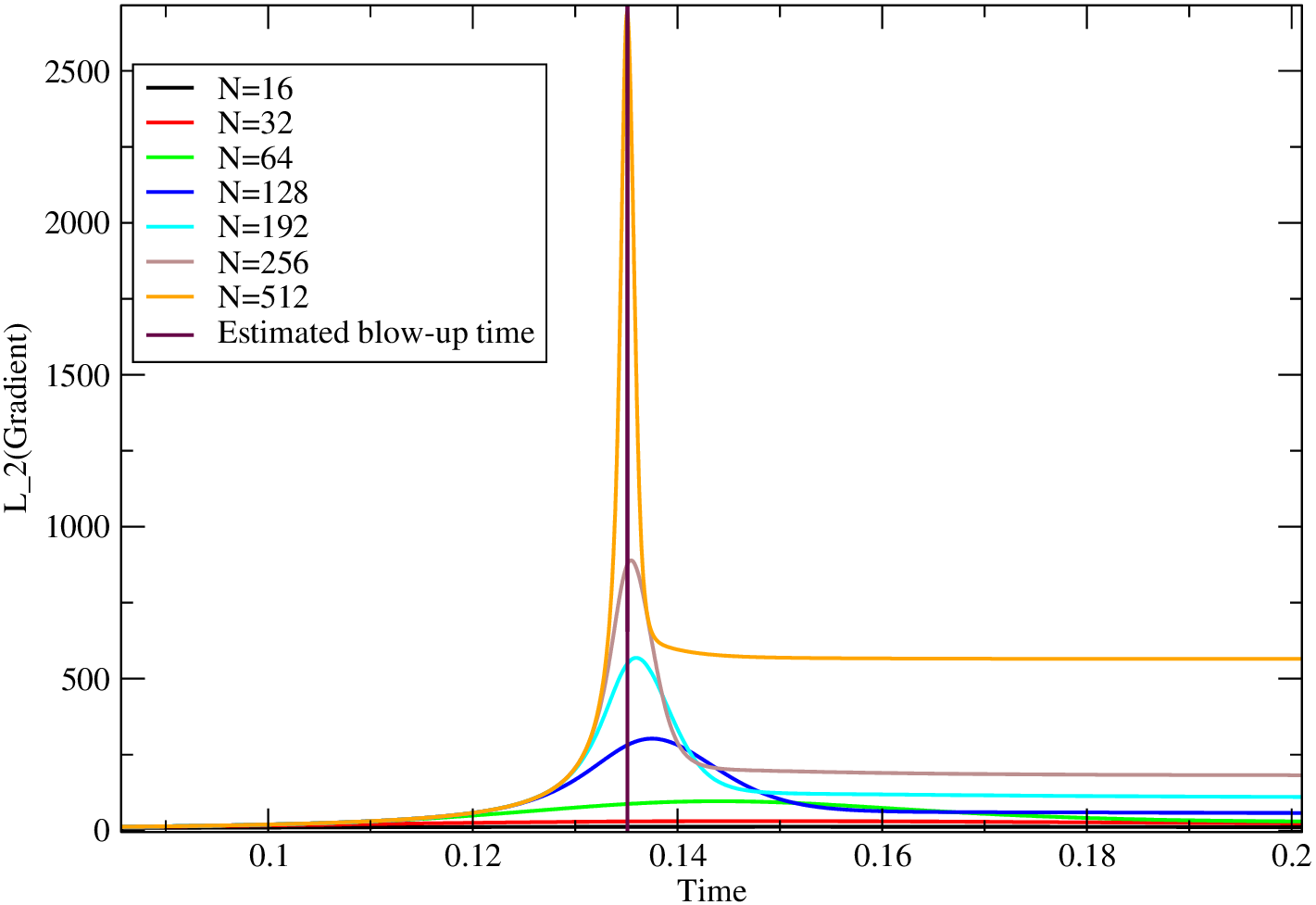}
\caption{Evolution of gradient $l_2$ norm for resolved modes. The vertical line denotes the numerically estimated blow-up instant $T=0.13504$ calculated with a mesh refinement algorithm.}
\label{critical_2}
\end{figure}

\begin{figure}
\centering
\includegraphics[width=3.5in]{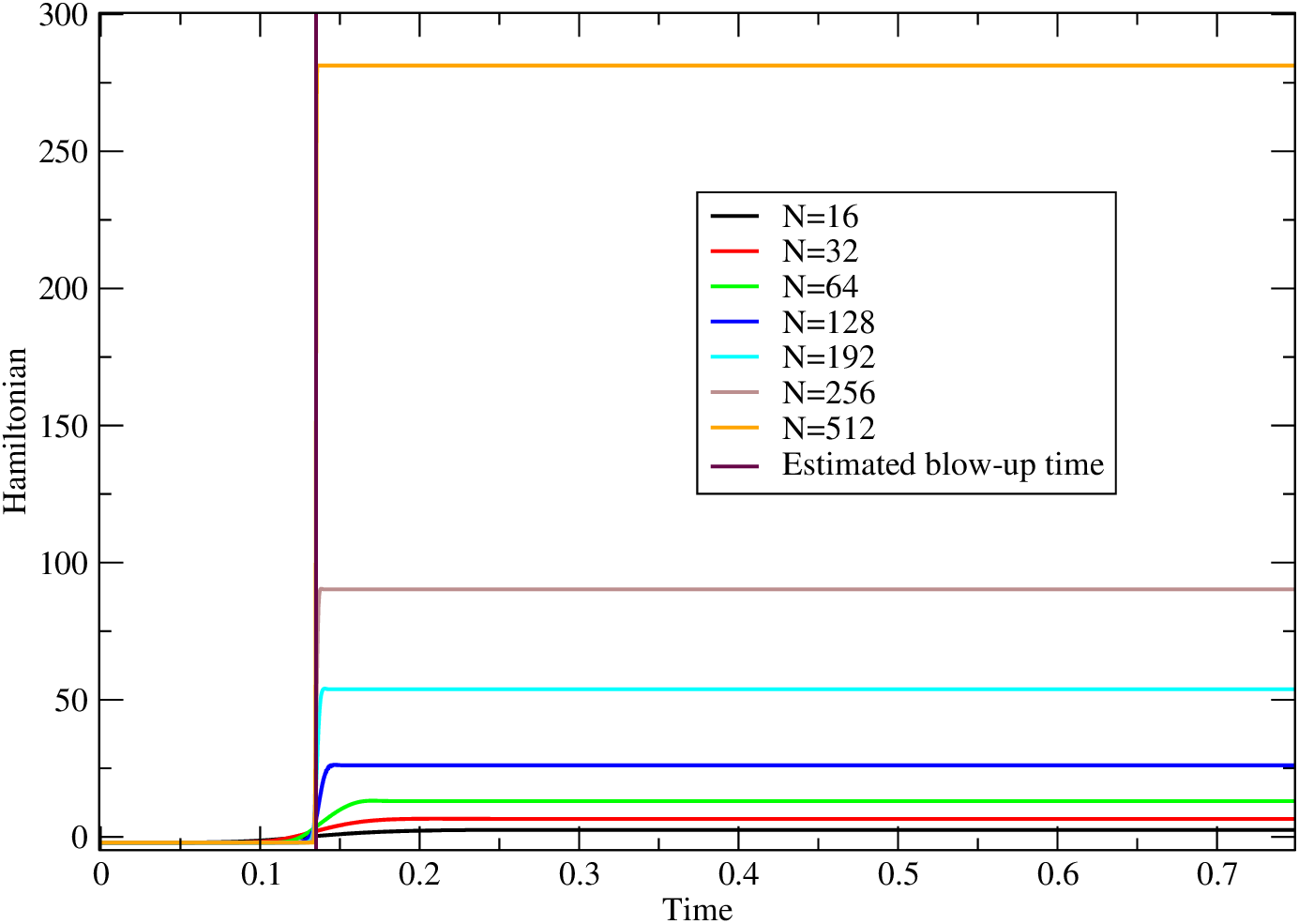}
\caption{Evolution of Hamiltonian for resolved modes. The vertical line denotes the numerically estimated blow-up instant $T=0.13504$ calculated with a mesh refinement algorithm.}
\label{critical_3}
\end{figure}

Figure \ref{critical_2} shows the evolution of the $l_2$ norm of the gradient for the resolved modes, i.e. $\sum_{k \in F} k^2 |u'_k|^2.$ As we have already said the $L_2$ norm of the gradient of the solution should blow up at the singularity. The behavior of the numerically computed $l_2$ norm of the gradient for the resolved modes is consistent with such a behavior and is again in remarkable agreement with the estimated time of the singularity. Similarly, as seen in Figure \ref{critical_3}, the value of the Hamiltonian has a jump of increasing magnitude with increasing resolution and the instant of the jump agrees again very well with the estimated blow-up instant.

\begin{figure}
\centering
\includegraphics[width=3.5in]{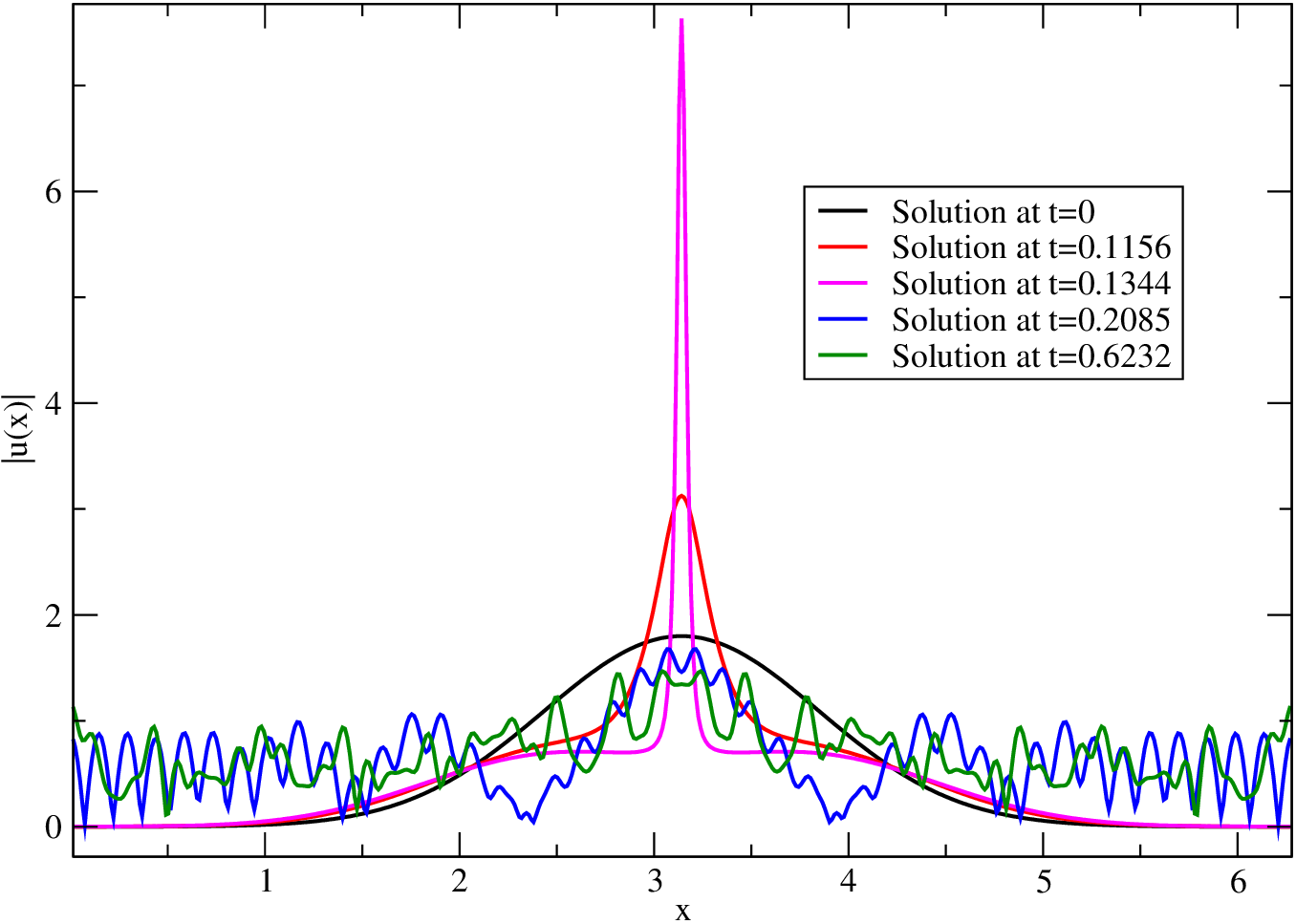}
\caption{Magnitude of the solution of the $t$-model with $N=512$ at different instants before and after the singularity (estimated singularity time is $T=0.13504$).}
\label{critical_4}
\end{figure}

\subsection{Smoothness of solution after the singularity}\label{smoothness}

Our purpose in this section is to provide more details which show that the solution computed from the $t$-model before, and especially after the singularity, is consistent with the theoretical results about the log-log blow up scenario in \cite{merle2} and the post-singular solutions in \cite{tao}. In particular, as it is shown in Merle and Raphael's work (see Theorem 3 in \cite{merle2}), for the log-log blow-up scenario, which our solution also follows, the quantization of mass at the singularity leads to a decomposition of the solution in a component that blows up in a self-similar fashion and another component, which is in $L_2$ but {\it not} in $H^1.$ In this case \cite{tao}, Tao's post-singular solutions exist but they are only in $L_2,$ i.e. they have lost regularity. Inevitably, in a numerical simulation we are bound to have finite resolution. Yet, we want to show that the trends exhibited by the $t$-model solution as we increase the resolution are consistent with the theoretical results. 

As we have seen in Figure \ref{critical_2}, the gradient $l_2$ norm seems to increase without bound (as it should) at the singularity as we increase the resolution. More importantly, its value after the singularity also increases with increased resolution. We have to note that the value of the gradient $l_2$ norm after the singularity is not a constant even though it may appear to be so from the figure. In fact, it oscillates mildly (about $0.04\%$) around a well-defined mean.  

Even more abrupt is the increase with resolution of the constant value of the Hamiltonian after the singularity (see Figure \ref{critical_3}). Since the mass value after the singularity appears to have converged with resolution (see Figure \ref{critical_1}), the abrupt increase of the Hamiltonian with resolution must be due to the increase of the gradient $l_2$ norm. If the trends shown in Figure \ref{critical_2} continue with even higher resolution, then the gradient $l_2$ norm should not only blow-up at the singularity but remain infinite after the singularity too. That would mean that the solution has lost its regularity.  

\begin{figure}
\centering
\includegraphics[width=3.5in]{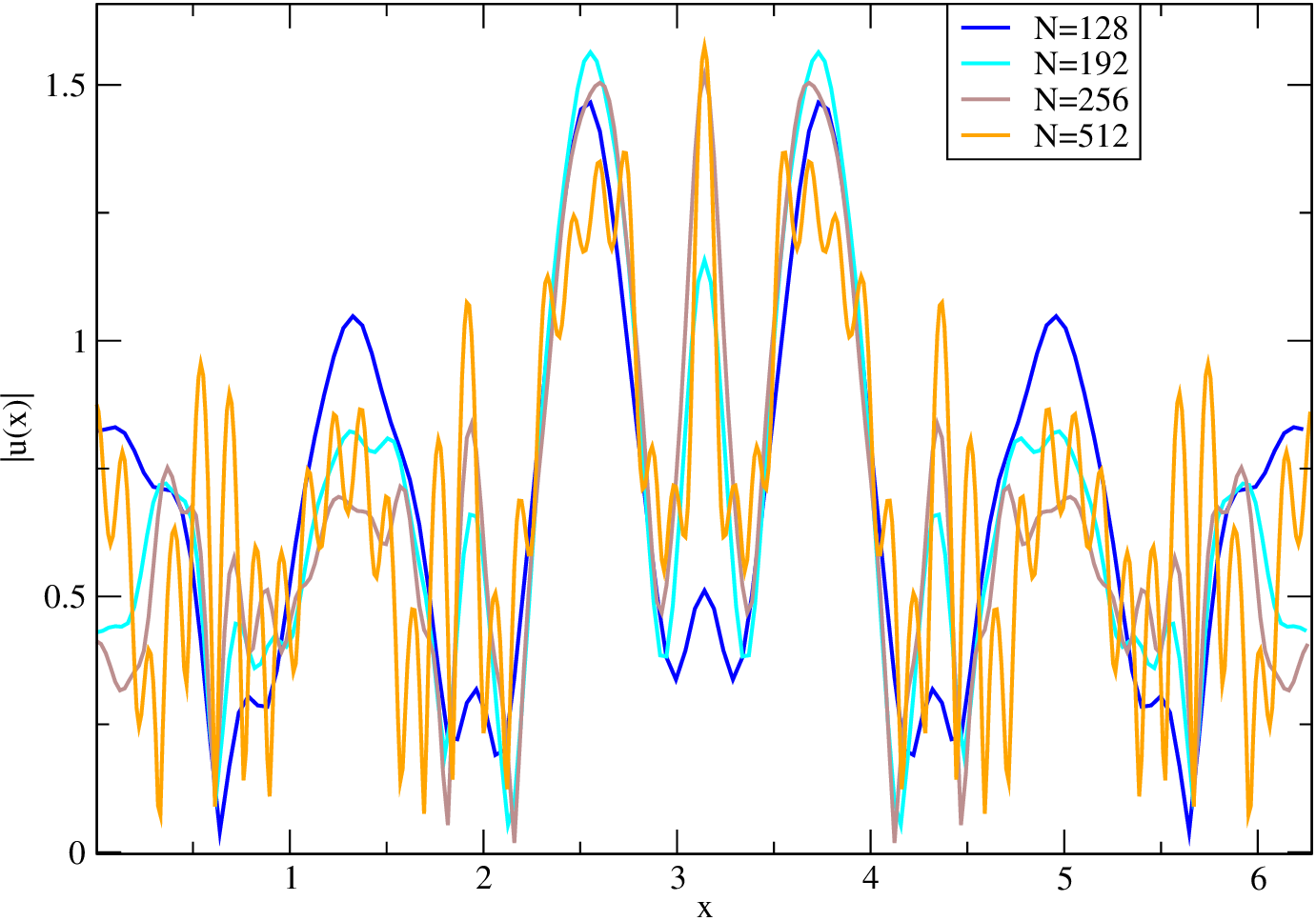}
\caption{Magnitude of the solution of the $t$-model with increasing resolution at time $t=0.75.$}
\label{critical_4_detail}
\end{figure}

Figure \ref{critical_4} shows the magnitude of the $t$-model solution with $N=512$ at different instants before and after the singularity. Even though the evolution leading to the singularity involves the concentration of the mass at one spatial point, long after the singularity has happened and mass has been ejected out of the resolved modes, the solution appears to be scattering towards the boundary of the spatial domain. We want to see how  the smoothness of the solution after the singularity changes as we increase the resolution. Figure \ref{critical_4_detail} shows the form of the solution for different resolutions at time $t=0.75.$ We see that as the resolution is increased the solution becomes rougher. This increase in roughness with resolution is consistent with the solution losing its regularity after the singularity. 

\begin{figure}
\centering
\includegraphics[width=3.5in]{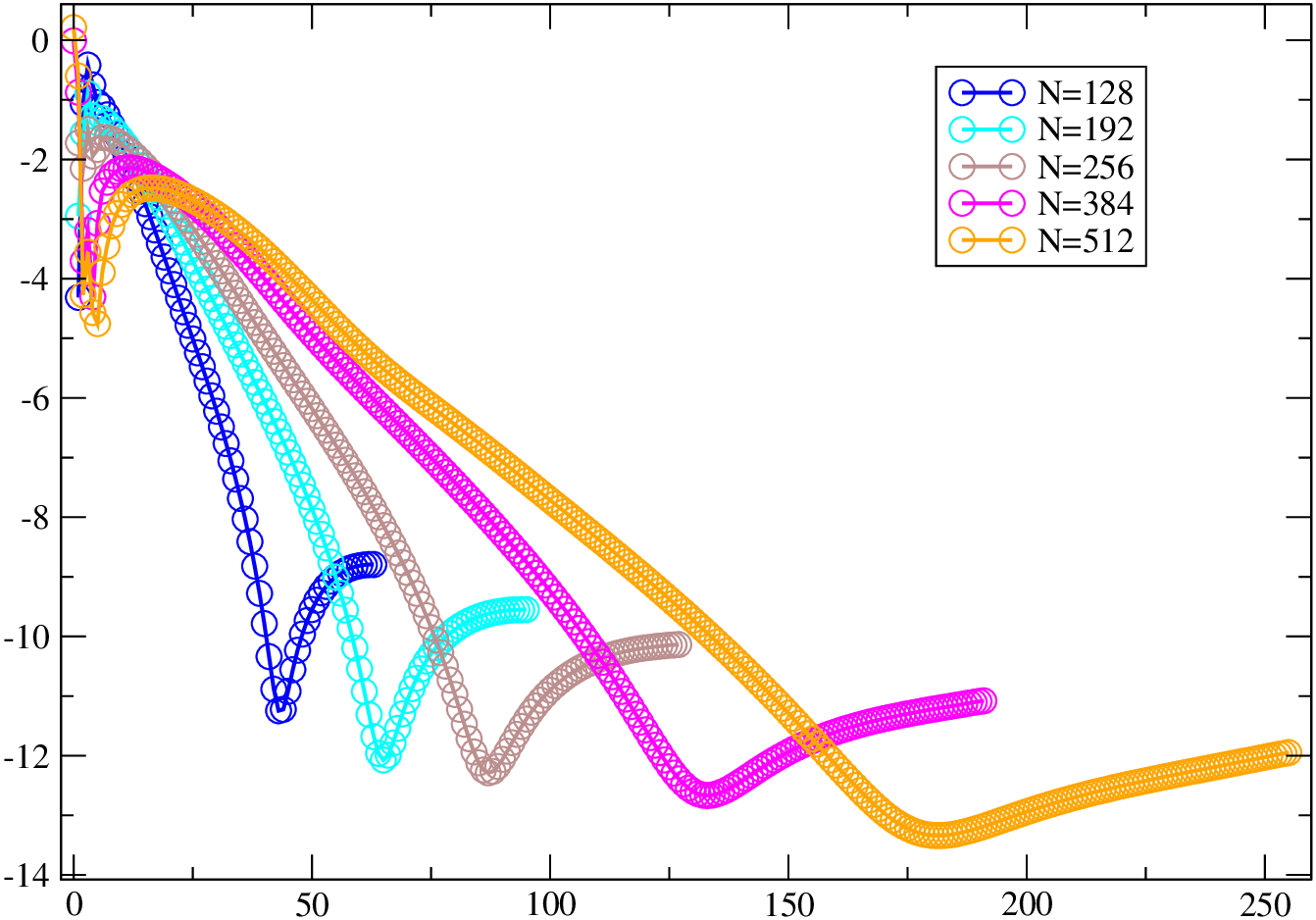}
\caption{Linear-log plot of the mass spectrum of the solution of the $t$-model with different resolutions at time $t=0.1355.$}
\label{critical_5}
\end{figure}

\begin{figure}
\centering
\includegraphics[width=3.5in]{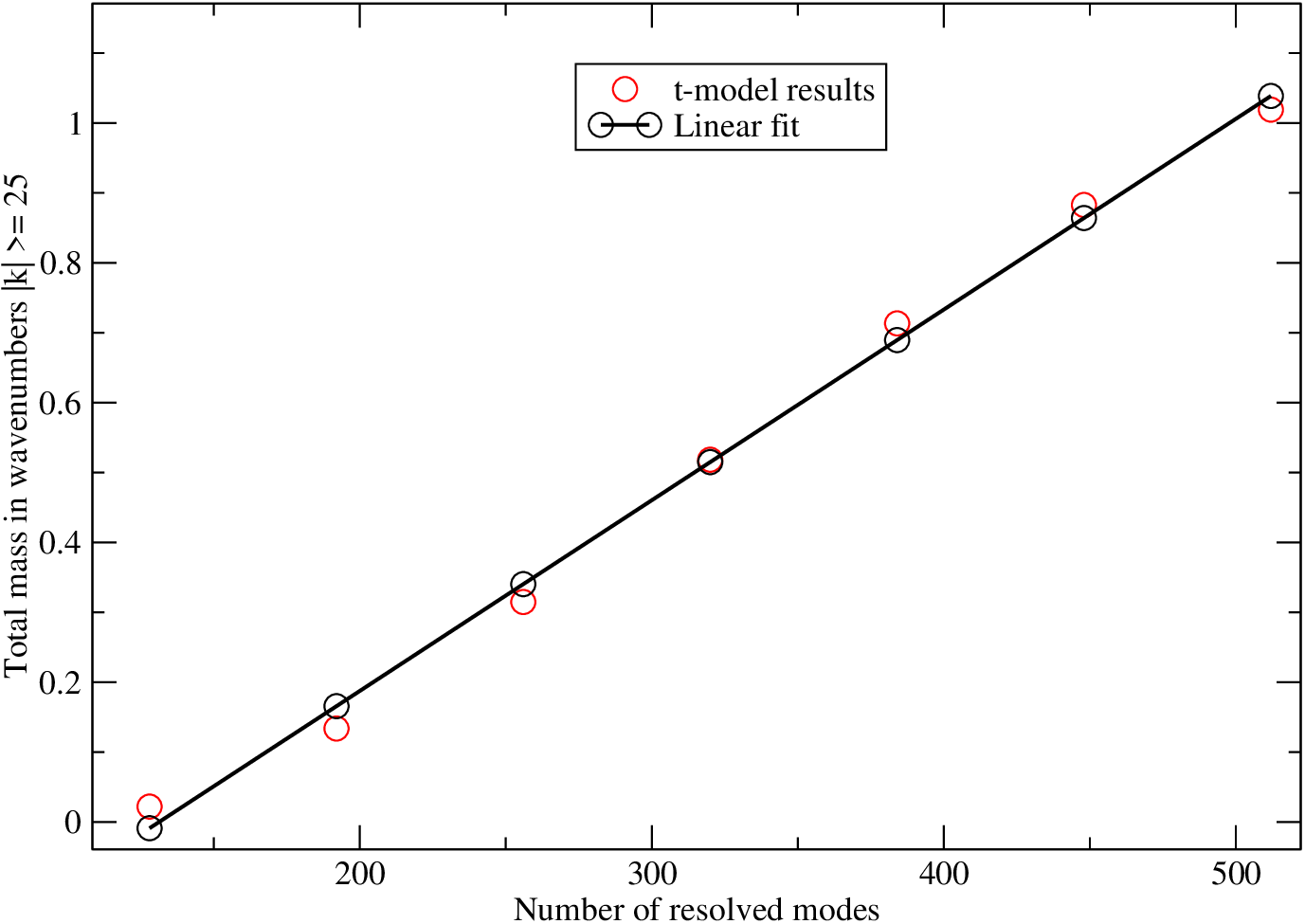}
\caption{Total mass in wavenumbers greater or equal to 25 as a function of the resolution at time $t=0.1355.$}
\label{critical_7}
\end{figure}

\begin{figure}
\centering
\includegraphics[width=3.5in]{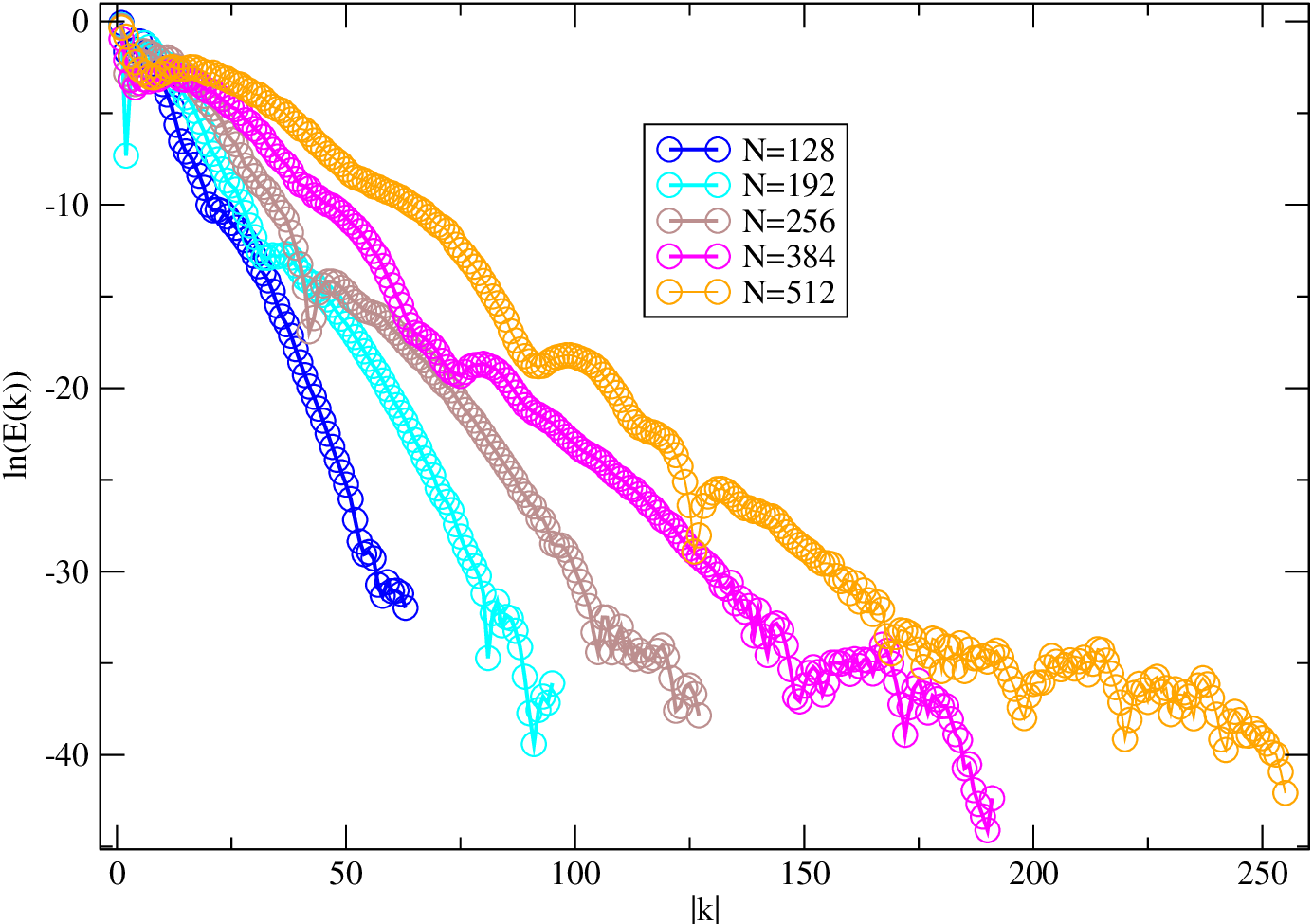}
\caption{Linear-log plot of the mass spectrum of the solution of the $t$-model with different resolutions at time $t=0.75.$}
\label{critical_6}
\end{figure}

To support our claim about the increase in roughness with increasing resolution we also look at the mass spectrum. In Figures \ref{critical_5} and \ref{critical_6} we show the spectrum of the $t$-model solution for different resolutions at time $t=0.1355$ and $t=0.75$ respectively. We make two observations. First, the $t$-model solution even after the singularity ($t=0.1355),$ remains well resolved. For the highest resolution here, $N=512,$ the 80 highest wavenumbers have a total mass of $O(10^{-5})$ while the total mass of the solution is $O(1).$ Second, the spectrum becomes shallower as we increase the resolution. 

The previous statement can be made more precise. In Figure \ref{critical_7} we have plotted the total mass of the resolved modes with wavenumber greater or equal to 25 shortly after the singularity ($t=0.1355$). The choice of the value 25 is because for the largest resolution here ($N=512$), this is the wavenumber where the linear tail (in linear-log coordinates) of the mass spectrum begins. We see that as the resolution is increased there is an increase in the total mass in the modes with wavenumber greater or equal to 25. The total mass in these modes as a function of the resolution can be fit by a straight line with a correlation coefficient of about 0.998. Since the total mass in all the modes appears to converge with increasing resolution (see Figure \ref{critical_1}), the increase in mass for the wavenumbers larger or equal to 25 must be attributed to a decrease in the mass of the wavenumbers less than 25. In other words, as the resolution is increased the large scales of the solution have a decreasing mass. Consequently, the solution becomes rougher.     

For $t=0.75,$ (see Figure \ref{critical_6}) the total mass of the 80 highest wavenumbers is $O(10^{-14})$ and the spectrum again becomes pronouncedly shallower with an increase in the resolution. A similar analysis as in the case shortly after the singularity shows that the total mass in the wavenumbers greater or equal to 25 also increases with increasing resolution. So, whether shortly or long after the singularity, as the resolution is increased the solution becomes rougher.

\section{Discussion}\label{discussion}
We have presented numerical results from the application of a reduced model, called the $t$-model, to the 1D critical Schr\"odinger equation for an initial condition that leads to a finite time singularity. The results suggest, in agreement with the recent theoretical results in \cite{tao}, that the singularity causes a finite amount of mass to be ejected instantaneously, after which, the solution continues to conserve the remaining mass. The time of occurrence of the singularity as predicted by the $t$-model is in very close agreement to the estimated singularity time given by a mesh refinement algorithm. It is very encouraging that the reduced model calculations are consistent with the blow-up instant predicted independently through a mesh refinement algorithm. Such consistency means that reduced models can help to shed light in open problems about the behavior of time dependent partial differential equations with singular (or near-singular) solutions.  

In addition to the instantaneous mass ejection, the $t$-model solution's roughness increases with increased resolution. This trend suggests that the solution loses regularity after the singularity. This is again in agreement with the theoretical results \cite{merle2}.   

Finally, similar behavior for the solution to the one shown in the current paper has been obtained for solutions of the critical Schr\"odinger equation in dimensions two and three. As we have already seen, even in one dimension, as we increase the resolution the calculation becomes rather expensive when the solution passes through the singularity. The calculations for the cases with more spatial dimensions are also very expensive as we increase the resolution and a detailed study definitely requires parallelization of the algorithm. This is ongoing work and the results will be presented elsewhere. 

\section*{Acknowledgements} I want to thank Prof. T. Tao for a helpful and enjoyable discussion during his visit at the University of Minnesota. I also want to thank Prof. V. Sverak for clarifying discussions about the behavior of solutions of nonlinear Schr\"odinger equations.

\end{document}